\documentclass[11pt,english]{article}
\usepackage[english]{babel}
\usepackage{amsfonts}
\usepackage{epsfig}
\usepackage{graphicx}
\usepackage{fancyhdr}
\usepackage{amssymb}
\usepackage{amsmath}
  \usepackage{verbatim} 
  \usepackage{appendix} 
\newtheorem{theorem}{Theorem}

\newtheorem{prop}{Proposition}[section]
\newtheorem{lemma}[prop]{Lemma}
\newtheorem{definition}[prop]{Definition}

\def\Ag {{\cal A}} 
\def\Lg {{\cal L}} 

\def\and {{\rm \; and \;}}

\def\exp {{\rm exp}}
\def\Ad{{A_1}} 
\def\sd{{s_{0,1}}} 
\def\Td{T_1}
\def\se{s_{0,2}}
\def\sf{s_{0,3}}
\def\Ag{A_4}
\def\sg{s_{0,4}}
\def\Ah{A_5}
\def\sh{s_{0,5}}
\def\si{s_{0,6}}
\def\Ti{T_6}
\def\Aj{A_7}
\def\sj{s_{0,7}}
\def\Tj{T_7}
\def\Ak{A_8}
\def\sk{s_{0,8}}
\def\Tl{T_9}

\newcommand {\R}{ \mathbb{R}}
\newcommand {\C}{ \mathbb{C}}
\newcommand {\N}{ \mathbb{N}}
\newcommand{\dd}{d_0, d_1, \tilde{d}_0,\tilde{d}_1}
\newcommand {\pa}{\partial}

\newcommand {\beqna} {\begin{eqnarray}}
\newcommand {\eeqna} {\end{eqnarray}}
\newcommand {\beqtn} {\begin{equation}}
\newcommand {\eeqtn} {\end{equation}}

\begin{document}

\title{\bf{Profile for a simultaneously blowing up solution to a complex valued  semilinear heat equation}}

\author{Nejla Nouaili,\\ Nejla.Nouaili@dauphine.fr\\ CEREMADE, Universit\'e Paris Dauphine, Paris Sciences et Lettres.\\ 
Hatem Zaag\footnote{This author is supported by the ERC Advanced Grant no. 291214,  BLOWDISOL and  by ANR Project ANA\'E ref. ANR-13-BS01-0010-03.},\\
Hatem.Zaag@univ-paris13.fr\\ Universit\'e Paris 13, Sorbonne Paris Cit\'e,\\ LAGA, CNRS (UMR 7539), F-93430, Villetaneuse, France.}

\maketitle
\begin{abstract}

We construct a solution to a complex nonlinear heat equation which blows up in finite time $T$ only at one blow-up point. We also give a sharp description of its blow-up profile. The proof relies on the reduction of the problem to a finite dimensional one  and the use of index theory to conclude. 
We note that the real and imaginary parts of the constructed solution blow up simultaneously, with the imaginary part dominated by the real.

\end{abstract}

\section{Introduction}
This paper is concerned with blow-up solutions of the complex heat equation
\beqtn
 \pa_t u=\Delta u+u^2,
\label{eqF}
\eeqtn
where $u(t):x\in \R^N\to\C$ and $\Delta$ denotes the Laplacian.\\
If we write $u(x,t)=v(x,t)+i\tilde{v}(x,t)$, where $v$ and $\tilde{v}\in \R$, we obtain the following reaction-diffusion system:
\beqtn
\begin{array}{lll}
\pa_t v&=&\Delta v+v^2-\tilde{v}^2, \\
\pa_t \tilde{v}&=&\Delta \tilde{v}+2 v \tilde{v},
\end{array}
\label{uv}
\eeqtn
 where $(x,t)\in\R^N\times (0,T) \mbox{, }v(0,x)=v_0(x)$ and $\tilde{v}(0,x)=\tilde{v}_0(x)$.\\
Equation (\ref{eqF}) has a strong relation with the viscous Constantin-Lax-Majda equation, which is a one dimensional model for the vorticity equation. For more details see Okamoto, Sakajo and Wunsch \cite{OSWN08}, Sakajo \cite{SJMSUT03} and \cite{SN03} and
Guo, Ninomiya, Shimojo and Yanagida in \cite{GHMY12}.

\medskip

The Cauchy problem for system (\ref{uv}) can be solved in $(L^\infty(\R^N))^2$, locally in time. We say that $u(t)=v(t)+i\tilde{v}(t)$ blows up in finite time $T<\infty$, if $u(t)$ exists for all $t\in [0,T)$ and 
$\lim_{t\to T}\|v(t)\|_{L^\infty}+\|\tilde{v}(t)\|_{L^\infty}=+\infty.$
In that case, $T$ is called the blow-up time of the solution. A point $x_0\in\R^N$ is said to be a blow-up point if there is a sequence $\{(x_j,t_j)\}$, such that $x_j\to x_0$, $t_j\to T$ and $|v(x_j,t_j)|+|\tilde{v}(x_j,t_j)|\to \infty$ as $j\to\infty$. The set of all blow-up points is called the blow-up set.

\bigskip

When $u$ is real (i.e., $\tilde{v}\equiv 0$), then system (\ref{uv}) is reduced to the scalar equation
\beqtn
 \pa_t u=\Delta u +u^p\mbox{, where } p=2.
\label{uscalar}
\eeqtn
The blow-up question for equation (\ref{uscalar}), with $p>1$, has been studied intensively by many authors and no list can be exhaustive.\\
When it comes to deriving the blow-up profile, the situation is completely understood in one space dimension (however, less is understood in higher dimensions, see Vel{\'a}zquez \cite{VCPDE92, VTAMS93,VINDIANA93} and Zaag \cite{ZIHP02, ZCMP02, ZMME02} for partial results). In one space dimension, given $a$ a blow-up point, this is the situation:

\medskip

{\textit {
\begin{itemize}
\item either
\beqtn
\sup_{|x-a|\leq K\sqrt{(T-t)\log(T-t)}}\left|(T-t)u(x,t)-f\left(\frac{x-a}{\sqrt{(T-t)\log(T-t)}}\right)\right|\to 0,
\label{profileu}
\eeqtn
\item or for some $m\in \N$, $m\ge 2$, and $C_m>0$
 \beqtn
\sup_{|x-a|<K(T-t)^{1/2m}}\left|(T-t)u(x,t)-f_m\left(\frac{C_m (x-a)}{(T-t)^{1/2m}}\right)\right|\to 0,
\label{um}
\eeqtn
as $t\to T$, for any $K>0$, where 
\begin{equation}
\label{deff}
f(z)=\frac 8{8+|z|^2}\mbox{ and }f_m(z) = \frac 1{1+|z|^{2m}} \mbox{, for all }z\in\R^N.
\end{equation}
\end{itemize}
}}

From Bricmont and Kupiainen \cite{BKN94} and Herrero and Vel{\'a}zquez \cite{HVAIHP93}, we have examples of initial data leading to each of the aboved-mentioned scenarios. Note that (\ref{profileu}) corresponds to the fundamental mode of the harmonic oscillator in the leading order, whereas (\ref{um}) corresponds to higher modes. Moreover, Herrero and 
Vel{\'a}zquez proved the genericity of the behavior (\ref{profileu}) in one space dimension in \cite{HVADSSP92} and \cite{HVCRAS94}, and only announced the result in the higher dimensional case (the result has never been published). Note also that the stability of such a profile with respect to initial data has been proved by Fermanian Kammerer, Merle and Zaag  in \cite{FZN00} and \cite{FMZMA00}. For more results on equation (\ref{uscalar}), see \cite{BQJMOS77}, \cite{GKCPAM85}, \cite{GKIUMJ87}, \cite{GKCPAM89}, \cite{HVAIHP93}, \cite{HVCRAS94}, \cite{MMCPAM04}, \cite{MMJFA08}, \cite{MZCPAM98}, \cite{MZMA00}, \cite{MMA07} and \cite{QSBV07}.

\bigskip

As we inferred above, equation (\ref{eqF}) appears as a complex generalization of the real valued equation (\ref{uscalar}). Note that there is another complex generalization of (\ref{uscalar}). Indeed, Filippas and Merle consider in \cite{FMJDE95} the following equation: 
\beqtn
\pa_t u =\Delta u+|u|^{p-1}u \mbox{ with }u\in\C\mbox{ and }p>1,
\label{eqFM}
\eeqtn
and generalize to this equation the results first proved in the real case by Giga and Kohn \cite{GKCPAM85,GKIUMJ87,GKCPAM89}. Our equation (\ref{eqF}) appears then as a ``twin`` of equation (\ref{eqFM}). However there is a fundamental difference between the two. Indeed, equation (\ref{eqFM}) has a variational structure, which allows to use various energy techniques, unlike equation (\ref{eqF}), where such techniques certainly fail. Therefore, considering equation (\ref{eqF}) appears as a highly challenging task.

\bigskip

Considering equation (\ref{eqF}) with non-real solutions, we have the following blow-up results from \cite{GHMY12}:

\medskip

\textit{(A) A non-simultaneous blow-up criterion, see Theorem 1.5 in \cite{GHMY12}:\\
 Assume that} 
\beqtn
v_0,\; \tilde{v}_0\in C^1(\R^N),\; 0\leq v_0\leq M,\; v_0\not\equiv M,\; 0<\tilde{v}_0\leq L, 
\label{CB1}
\eeqtn

\beqtn
\lim_{|x|\to\infty} v_0(x)=M,\;\;\lim_{|x|\to\infty }\tilde{v}_0(x)=0,
\label{CB2}
\eeqtn 
\textit{for some constants $L>0$ and $M > 0$. Then, the solution of (\ref{uv}) blows up at time $t=T(M)$ with $\tilde{v}\not\equiv 0$. Moreover, the component $v$ blows up only at space infinity and $\tilde{v}$ remains bounded.
}

\medskip
\textit{(B) A Fourier-based blow-up criterion, see Theorem 1.2 in \cite{GHMY12}:\\
If the Fourier transform of initial data of (\ref{eqF}) is real and positive, then the solution blows up.}

\medskip

\textit{(C) A simultaneous blow-up criterion, see  Theorem 1.3 in \cite{GHMY12}:\\
If $v_0$ is even, $\tilde{v}_0$ is odd with $\tilde{v}_0(x)>0$ for $x>0$, then the fact that the blow-up set is compact implies that $v$ and $\tilde{v}$ blow up simultaneously.}

\bigskip

Following the description in (\ref{profileu}) and (\ref{um}) together with the work of \cite{GHMY12}, we see that the blow-up profile derivation remains open, in the non-real case. Indeed, either we have the description 
(\ref{profileu}) or (\ref{um}), with a zero imaginary part, or we have blow-up solutions from  \cite{GHMY12}, with a non-trivial imaginary part, and no profile description.

\medskip

In this paper, we give the first example of a complex-valued blow-up solution of equation \eqref{eqF},  with a non-trivial imaginary part, and a full description of its blow-up profile, obeying behavior \eqref{profileu} (note that our method extends with no difficulty to the construction of an analogous solution obeying the behavior \eqref{um}; however, the proof should be even more technical). Let us note that the blow-up behavior we give here is not predicted by \cite{GHMY12} (see details in the remarks following our result). More precisely, this is our result:
\begin{theorem} ({\bf{Existence of a blow-up solution for equation (\ref{eqF}) with the description of its profile}}) There exists $T>0$ such that equation (\ref{eqF}) has a solution $u(x,t)=v(x,t)+i\tilde{v}(x,t)$ in $\R^N\times [0,T)$ such that:\\
(i) the solution $u$ blows up in finite time $T$ only at the origin; \\
(ii) there holds that
\beqtn
\left \| (T-t)u(.,t)-f \left(\frac{.}{\sqrt{(T-t)|\log(T-t)|}}\right)\right \|_{L^{\infty}}\leq \frac{C}{\sqrt{|\log(T-t)|}},
\label{profilev}
\eeqtn 
where $f$ is defined by (\ref{deff}).\\
(iii) For all $R>0$,
\beqtn
\sup_{|x|\leq R\sqrt{T-t}}\left| (T-t)\tilde{v}(x,t)-\frac{\sum_{i=1}^{N}C_i}{|\log(T-t)|^2}\left(\frac{x_{i}^{2}}{T-t}-2\right)\right|\leq \frac{C}{|\log(T-t)|^\alpha},
\label{profilevt}
\eeqtn
where $(C_1,C_2,..,C_N)\neq (0,0,..,0)$ and $2<\alpha\leq 2+\eta$ for some small $\eta>0$.\\
(iv) For all $x\neq 0$, $u(x,t)\to u^*(x)$ uniformly on compacts sets of $\R^N\backslash \{0\}$, and

\beqtn
u^*(x)\sim \frac{16 | \log|x||}{|x|^2} \mbox{ as }x\to 0.
\label{u*}
\eeqtn

\label{theorem1}
\end{theorem}
{\textbf{Remarks:}}

\medskip

1) Note that the real and imaginary parts of $u$ blow up simultaneously at $x=0$. However the real part dominates the imaginary part in the sense that 
\[v(0,t)\sim \frac{1}{T-t}>>\frac{2\left |\sum_{i=1}^{N}C_i\right|}{(T-t)|\log(T-t)|^2}\sim |\tilde{v}(0,t)|\mbox{ as }t\to\infty.\]

2) As announced right before the statement of our theorem, the solution we construct is new and doesn't obey the criteria given in \cite{GHMY12}: this is clear from (\ref{initialq}) below.
The proof relies on the reduction of the problem to a $2(N+1)-$dimensional problem (a $4-$dimensional one if $N=1$; see Lemma \ref{transversality} below).
In the real case treated by Merle and Zaag in \cite{MZDuke97}, the problem was of dimension $N+1$. Since that number is equal to the dimension of the blow-up parameters ($1$ for the blow-up time and $N$ for the blow-up point), the authors of \cite{MZDuke97} were able to show the stability of the behavior (\ref{profilev}) with respect to initial data, of course in the real case. Here, in the complex case, since the dimension of our problem $2(N+1)$ exceeds that of the blow-up parameters $(N+1)$, we suspect our solution to be unstable with respect to perturbations in initial data.

\bigskip

Our proof uses some ideas developed by  Bricmont and Kupiainen \cite{BKN94} and Merle and Zaag \cite{MZDuke97} to construct a blow-up solution for the semilinear heat equation (\ref{uscalar}) obeying the behavior (\ref{profileu}). In \cite{EZSMJ11}, Ebde and Zaag use the same ideas to show the persistence of the profile (\ref{profileu}) under perturbations of equation (\ref{eqF}) in the real case by lower order terms involving $u$ and $\nabla u$.
 In \cite{MZ07}, Masmoudi and Zaag adapted that method to the case of the following complex Ginzburg-Landau equation, where no gradient structure exists:
\[\pa_t u=(1+i\beta)\Delta u+(1+i\delta)|u|^{p-1}u\mbox{, with $\beta$ and $\delta$ are reals},\]
(note that the case $\beta=0$ and $\delta$ small was first studied by Zaag in \cite{ZAIHPANL98}).\\
More precisely, the proof relies on the understanding of the dynamics of the selfsimilar version of (\ref{uv}) (see system (\ref{qqtilde}) below) around the profile (\ref{profileu}). Moreover, we proceed in two steps:
\begin{itemize}
\item First, we reduce the question to a finite-dimensional problem: we show that it is enough to control a $(N+1)$-dimensional variable in order to control the solution (which is infinite dimensional) near the profile.
\item Second, we proceed by contradiction to solve the finite-dimensional problem and conclude using index theory.
\end{itemize}
Surprisingly enough, we would like to mention that this kind of methods has proved to be successful in various situations including hyperbolic and parabolic PDE, in particular with energy-critical exponents. This was the case for the construction of multi-solitons for the semilinear wave equation in one space dimension by C\^{o}te and Zaag \cite{CZCPAM13}, the wave maps by Rapha\"{e}l and Rodnianski \cite{RRIHES12}, the Schr\"{o}dinger maps by Merle, Rapha\"{e}l and Rodnianski \cite{MRRCRAS11}, the critical harmonic heat flow by Schweyer \cite{SJFA12} and the two-dimensional Keller-Segel equation by Rapha\"{e}l and Schweyer \cite{RSCPAM13}.

\medskip

We proceed in 3 sections to prove Theorem \ref{theorem1}. We first give in Section 2 an equivalent formulation of the problem in the scale of the well-known similarity variables. Section 3 is devoted to the proof of the similary variables formulation (this is a central part in our argument). Finally, we conclude the proof of Theorem \ref{theorem1} in Section 4.\\
{\bf Acknowledgement}. We would like to thank the referee for his valuable remarks which undoubtedly improved the presentation of our paper. 

 \section{Formulation of the problem}
For simplicity, we give the proof in one dimension. The adaptation to higher dimensions is straightforward. We would like to find initial data  $u_0=v_0+i\tilde{v}_0$ such that the solution $u=v+i\tilde{v}$ of equation (\ref{uv}) blows up in time $T$ with

\beqtn
\displaystyle\lim_{t\to T}\left\|(T-t)u(x,t)-f\left(\frac{x}{\sqrt{(T-t)|\log(T-t)|}}\right) \right\|_{L^\infty}=0,
\label{Ineuf}
\eeqtn
where $f$ is defined in (\ref{deff}).

\bigskip

This is the main estimate and the other results of Theorem \ref{theorem1} will appear as by-products of the proof (see Section 4 for the proof of all the estimates of Theorem \ref{theorem1}).\\
Introducing the following self-similar transformation of problem (\ref{uv}):
\beqtn
\begin{array}{c}
w(y,s)=(T-t)v(x,t)\mbox{, }\tilde{w}(y,s)=(T-t)\tilde{v}(x,t),\\
y=\frac{x-a}{\sqrt{T-t}}\mbox{, }s=-\log(T-t),
\end{array}
\label{chauto}
\eeqtn

we see that (\ref{Ineuf}) is equivalent to finding $s_0>0$ and initial data at $s_0$, $W_0(y,s_0)=w_0(y,s_0)+i\tilde{w}_0(y,s_0)$, such that $W(y,s)=w(y,s)+i\tilde{w}(y,s)$
satisfies 
\beqtn
\displaystyle\lim_{s\to\infty}\left\|W(y,s)-f\left(\frac{y}{\sqrt{s}}\right)\right\|_{L^\infty}=0.
\label{ProfileW}
\eeqtn
Introducing 
\beqtn
w=\varphi+q\mbox{ and }\tilde{w}=\tilde{q}\mbox{ where }\varphi(y,s)=f\left(\frac{y}{\sqrt{s}}\right)+\frac{1}{4 s},
\label{initialqqtilde}
\eeqtn
the problem is then reduced to constructing a function $Q=q+i\tilde{q}$ such that
\[\lim_{s\to\infty}\|Q(y,s)\|_{L^\infty}=0,\]
and $(q,\tilde{q})$ is a solution of the following equation for all $(y,s)\in \R\times[s_0,\infty)$,
\beqtn
\begin{array}{lll}
\displaystyle \pa_s q&=&\displaystyle({\cal{L}}+V)q+b(y,s)+R(y,s),\\
\displaystyle \pa_s \tilde{q}&=&\displaystyle({\cal{L}}+V)\tilde{q}+\tilde{b}(y,s),
\end{array}
\label{qqtilde}
\eeqtn
where
\beqtn
\displaystyle{\cal{L}}=\pa_{y}^{2} -\frac{1}{2}y\pa_y+1\mbox{, }V(y,s)=2\left(\varphi(y,s)-1\right),
\label{OperatorL}
\eeqtn
\beqtn
b(y,s)=q^2-\tilde{q}^2\mbox{, }\tilde{b}(y,s)=2q\tilde{q},\label{defN}
\eeqtn
and
\beqtn
\displaystyle R(y,s)=\pa_{y}^{2} \varphi-\frac{1}{2}y \pa_y \varphi-\varphi+\varphi^2-\pa_s \varphi.
\eeqtn

The control of $q$ and $\tilde{q}$ near $0$ obeys two facts:
\begin{itemize}
\item Localization: the fact that our profile $\varphi(y,s)$ dramatically changes its value from $1+\frac{1}{4s}$ in the region near $0$ to $\frac{1}{4s}$ in the region near infinity, according to a free boundary moving at the speed $\sqrt{s}$. This will require different treatments in the regions $|y|<2K_0\sqrt{s}$ and $|y|>2K_0\sqrt{s}$ for some $K_0$ to be chosen.

\item Spectral information: the fact that the operator $\Lg$ is selfadjoint, $b$ and $\tilde{b}$ are quadratic in $(q,\tilde{q})$ and  that
\[\|R(s)\|_{L^\infty}+\|V(s)\|_{L^{2}_{\rho}}\to 0 \mbox{ as } s\to \infty
\]
from (\ref{initialqqtilde}) and (\ref{OperatorL}), which shows that the dynamics of equation (\ref{qqtilde}) near $0$ are driven by the spectral properties of $\Lg$. This will require a decomposition of the solution according to the spectrum of $\Lg$. Note that the operator $\Lg$ is self-adjoint in the Hilbert space
\[
L^{2}_{\rho}=\lbrace g\in L^{2}_{loc}(\R,\C)\mbox{,
}\|g\|_{L^{2}_{\rho}}^{2}\equiv\int_{\R}|g|^2e^{-\frac{|y|^2}{4}}dy
<+\infty\rbrace\mbox{ where }\rho(y)=\displaystyle\frac{e^{-\frac{|y|^2}{4}}}{(4\pi)^{1/2}}.
\]
\end{itemize}
The spectrum of $\Lg$ is explicitly given by 
\[spec(\Lg)=\{1-\frac{m}{2}\mbox{, }m\in\N\}.\]
All the eigenvalues are simple, the eigenfunctions are dilations of Hermite's polynomial and  given by 
\beqtn
h_m(y)=\sum_{n=0}^{[\frac{m}{2}]}\frac{m!}{n!(m-2n)!}(-1)^n y^{m-2n}.
\label{hermite}
\eeqtn 

Note that $\Lg$ has two positive (or expanding) directions ($\lambda=1$ and $\lambda=\frac 1 2$), and a zero direction ($\lambda=0$).
Complying with the localization and spectral information facts, we will decompose $q$ and $\tilde{q}$ accordingly as stated above:
\begin{itemize}
\item first, we consider a non-increasing cut-off function $\chi_0 \in C^{\infty}_{0}(\R^+,[0,1])$ such that $supp(\chi_0)\subset [0,2]$, $\chi_0(\xi)=1$ for $\xi<1$ and $\chi_0(\xi)=0$ for $\xi>2$, then introduce
\[\chi(y,s)=\chi_0\left(\frac{|y|}{K_0\sqrt{s}}\right),\]
where $K_0\geq 1$ will be chosen large enough so that various technical estimates hold. Then, we write $q=q_b+q_e$ and $\tilde{q}=\tilde{q}_b+\tilde{q}_e$, where the inner parts and the outer parts are given by 
\[q_b=q\chi\mbox{, }\tilde{q}_b=\tilde{q}\chi\mbox{, } q_e=q(1-\chi)   \mbox{ and } \tilde{q}_e=\tilde{q}(1-\chi)  .\]
Let us remark that 
\[supp(q_b(s))\subset B(0,2K_0\sqrt{s})\mbox{ and }supp(q_e(s))\subset \R\setminus B(0,K_0\sqrt{s}),\]
and the same holds for $\tilde{q}_b$ and $\tilde{q}_e$.
\item Second, we study $q_b$ and $\tilde{q}_b$ using the structure of $\Lg$, isolating the nonnegative directions. More precisely we decompose $q_b$ and $\tilde{q}_b$ as follows
 \beqtn
 \begin{array}{lll}
 q_b(y,s)&=&\sum_{0}^{2}q_m(s) h_m(y)+q_-(y,s),\\
 \tilde{q}_b(y,s)&=&\sum_{0}^{2}\tilde{q}_m(s) h_m(y)+\tilde{q}_-(y,s),\\
 \label{decompositionq}
 \end{array}
 \eeqtn
 where $q_m$ (respectively $\tilde{q}_m$) is the projection of $q_b$ (respectively $\tilde{q}_b$) on $h_m$, $q_-(y,s)=P_-(q_b)$ (respectively $\tilde{q}_-(y,s)=P_-(\tilde{q}_b)$) and $P_-$ is the projection on $\{h_i,\; i\geq 3\}$ the negative subspace of the $\Lg$. 
\end{itemize}
In summary, we can decompose $q$ (respectively $\tilde{q}$) in 5 components as follows:
 \beqtn
 \begin{array}{lll}
 q(y,s)&=& \sum_{m=0}^{2}q_m(s) h_m(y)+ q_-(y,s)+q_e(y,s),\\
  \tilde{q}(y,s)&=& \sum_{m=0}^{2}  \tilde{q}_m(s) h_m(y)+   \tilde{q}_-(y,s)+  \tilde{q}_e(y,s).
 \end{array}
 \label{decompq}
 \eeqtn
 Here and throughout the paper, we call  $q_-(y,s)$ (respectively $\tilde{q}_-$) the negative part of $q$ (respectively $\tilde{q}$),  $q_2$ (respectively $\tilde{q}_2$), the null mode of $q$ (respectively $\tilde{q}$). 

 \label{Section2}
 \section{The construction method in selfsimilar variables}
 This section is devoted to the proof of the existence of a solution $(q,\tilde{q})$ of system (\ref{qqtilde}) satisfying $\|q(s)\|_{L^\infty}+\|\tilde{q}(s)\|_{L^\infty}\to 0$. This is a central argument in our proof. In Section 4, we use this solution and give the proof of Theorem \ref{theorem1}. Though we refer to the earlier work by Merle and Zaag \cite{MZDuke97} for purely technical details, we insist on the fact that we can completely split from that paper as long as ideas and arguments are considered. We hope that the explanation of the strategy we give in this section will be more reader friendly.

\bigskip

We proceed in 3 subsections:\\
%
%
 - In the first subsection, we give all the arguments of the proof without the details, which are left for the following subsection (readers not interested in technical details may stop here).\\
- In the second subsection, we give various estimates concerning initial data.\\
- In the third subsection, we give the dynamics of system (\ref{qqtilde}) near the zero solution, in accordance with the decomposition (\ref{decompq}).

\subsection{The proof without technical details}

Given $s_0>0$, we consider initial data for equation (\ref{qqtilde}) of the following form:
\beqtn
\begin{array}{lll}
\displaystyle q_{d_0,d_1}(y,s_0)&=&\displaystyle \frac{A}{s_{0}^{2}} (d_0+d_1 y)\chi(2y,s_0),\\
\displaystyle \tilde{q}_{\tilde{d}_0,\tilde{d}_1,\tilde{d}_2}(y,s_0)&=&\displaystyle \left[ \frac{\tilde{A}}{s_{0}^{\alpha}} (\tilde{d}_0+\tilde{d}_1 y)+\frac{\tilde{B}}{s_{0}^{2}}h_2(y)\right]\chi(2y,s_0),
\end{array}
\label{initialq}
\eeqtn
for some constants  $A$, $\tilde{A}$ and $\tilde{B}$ will be fixed later and the parameters $(\dd)\in [-2,2]^4$.
The solution of equation (\ref{qqtilde}) with initial data (\ref{initialq}) will be denoted by\\ $(q,\tilde{q})(s_0,A,\tilde{A},\tilde{B},(\dd),y,s)$, or, when there is no ambiguity by\\ $(q,\tilde{q})(s_0,\dd,y,s)$ or even  $(q,\tilde{q})(y,s)$. We will show that given $\tilde{B}\in\R$, if $A$ and $\tilde{A}$ are fixed large enough, then $s_0$ is fixed large enough depending on $A$, $\tilde{A}$ and $\tilde{B}$, we can also fixe the parameters $(\dd) \in[-2,2]^4$, so that the solution\\ $(q,\tilde{q})\left(s_0,A\,\tilde{A},\tilde{B},\dd,y,s\right)$ will converge to $0$ as $s\to\infty$.
Thanks to the decomposition given in (\ref{decompq}), in order to control $(q,\tilde{q})(s)$ near $(0,0)$, it is enough to control it in some shrinking set defined as follows:

\medskip

\begin{definition}\label{defVA}{\bf (Definition of a shrinking set for the components of $(q,\tilde{q})$)}
\noindent For all $A\geq 1$,  $\tilde{A}\geq 1$, $0<\eta<\frac 1 2$, $2<\alpha<2+\eta$ and $s\geq e$, we define $V_A(s)$  (respectively $\tilde{V}_{\tilde{A}}(s)$) 
  as the set of all functions $r$ (respectively $\tilde{r}$) in $L^\infty$ such that:
 \[|r_m(s)|\leq A s^{-2}\; m=0,1,\;\;  |r_2(s)|\leq A^2 (\log s) s^{-2},\]
 \[\forall y\in\R,\; |r_-(y,s)|\leq A(1+|y|^3)s^{-2},\; \|r_e(s)\|_{L^\infty}\leq A^2 s^{-\frac{1}{2}}, \] 
(respectively 
 \[|\tilde{r}_m(s)|\leq \tilde{A} s^{-\alpha}\; m=0,1,\;\;  |\tilde{r}_2(s)|\leq \tilde{A}^2  s^{-2+\eta},\]
 \[\forall y\in\R,\; |\tilde{r}_-(y,s)|\leq \tilde{A}(1+|y|^3)s^{-\alpha},\; \|\tilde{r}_e(s)\|_{L^\infty}\leq \tilde{A}^2 s^{-\alpha+3/2}), \] 
where $r_-$, $r_e$ and $r_m$ (respectively  $\tilde{r}_-$, $\tilde{r}_e$ and $\tilde{r}_m$) are defined in Section \ref{Section2}.
\end{definition}

\noindent As a matter of fact, if $s\geq e$ and $(r,\tilde{r})\in V_{A}(s)\times \tilde{V}_{\tilde{A}}(s)$, then one easily derives that

\beqtn
\|r(s)\|_{L^\infty}\leq C\frac{A^2}{\sqrt{s}}\mbox{ and }\|\tilde{r}(s)\|_{L^\infty}\leq C\frac{\tilde{A}^2}{s^{\alpha-3/2}},\label{petit}
\eeqtn
(see Proposition 3.7 page 157 in \cite{MZDuke97} for details). Thus, our aim become the following: 
\begin{prop}\label{Newprop}
{\bf (Existence of a solution of (\ref{qqtilde}) trapped in $V_A(s)\times\tilde{V}_{\tilde{A}}(s)$)}    
There exists $\Ad$ such that for all $A\geq \Ad$ and $\tilde{A}\geq \Ad$, $0<\eta<\frac{1}{10}$ and $2<\alpha<2+\eta$, there exists 
$\sd (A, \tilde{A})$ 
such that for all $s_0\geq \sd$ and $|\tilde{B}|\leq 1$, there exists $(d_0,d_1,\tilde{d}_0,\tilde{d}_1)$, such that,\\
if $(q, \tilde{q})$ is a solution of (\ref{qqtilde}) with initial data at $s_0$ given by (\ref{initialq}), then 
\[\forall s\geq s_0,\;\; q(s)\in V_A(s)\mbox{ and }\tilde{q}(s)\in \tilde{V}_{\tilde{A}}(s).\] 
\end{prop}
The aim of this section is to prove this proposition.

\bigskip

 In the following lemma, we find a set $D_{A,\tilde{A},\tilde{B},\eta,\alpha,s_0}=D_{s_0}$ such that $(q,\tilde{q})(s_0)\in  V_{A}(s_0)\times \tilde{V}_{\tilde{A}}(s_0)$, whenever $(\dd)\in D_{s_0}$. More precisely, we claim the following:

\begin{lemma}{\bf (Choice of parameters $\dd$ to have initial data in $V_A(s)\times\tilde{V}_{\tilde{A}}(s)$ at $s=s_0$)}  \label{initialisationN} For each $|\tilde{B}|\leq 1$, $A\geq 1$, $\tilde{A}\geq 1$, $0<\eta<\frac{1}{10}$ and $2<\alpha<2+\eta$, there exists $\se(A,\tilde{A},\tilde{B})\geq e$ 
such that for all $s_0\geq \se$:\\
If initial data for equation (\ref{qqtilde}) are given by (\ref{initialq}):
 then, there exists a cuboid
\beqtn
D_{A,\tilde{A},\tilde{B},\eta,\alpha,s_0}=D_{s_0} \subset [-2,2]^4,
\eeqtn
 such that, for all $(d_0,d_1,\tilde{d}_0,\tilde{d}_1)\in D_{s_0}$, we have 
\[(q,\tilde{q})(s_0,A,\tilde{A},\tilde{B},\dd)\in V_A(s_0)\times \tilde{V}_{\tilde{A}}(s_0).\]
\end{lemma}
\textit{Proof}: The proof is purely technical and follows as the analogous step in \cite{MZDuke97}, for that reason we refer the reader to Lemma 3.5 page 156 and Lemma 3.9 page 160 in \cite{MZDuke97}. $\blacksquare$

\bigskip

Let us consider $|\tilde{B}|\leq 1$, $A\geq 1$, $\tilde{A}\geq 1$, $0<\eta<1/10$, $2<\alpha<2+\eta$, $(\dd)\in D_{s_0}$ and $s_0\geq s_{0,1}$ defined in Lemma  \ref{initialisationN}. From the local Cauchy theory, we define a maximal solution $(q,\tilde{q})$ to equation (\ref{qqtilde}) with initial data (\ref{initialq}), and a maximal time $s_*(\dd)\in [s_0,+\infty]$ such that, for all $s\in[s_0,s_*)$, $(q,\tilde{q})(s)\in V_{A}(s)\times \tilde{V}_{\tilde{A}}(s)$ and:
\begin{itemize}
\item either $s_*=\infty$,
\item or $s_*<\infty$ and from continuity, $(q,\tilde{q})(s_*)\in \pa\left(V_{A}(s_*)\times \tilde{V}_{\tilde{A}}(s_*)\right)$, in the sense that when $s=s^*$, one '$\leq$' symbol in the definition of $V_{A}(s_*)$ and $\tilde{V}_{\tilde{A}}(s_*)$ is replaced by the symbol '$=$'. 
\end{itemize}
Our aim is to show that for all $|\tilde{B}|\leq 1$, for $A$, $\tilde{A}$ and $s_0$ large enough, one can find a parameter $(\dd)$ in $D_{s_0}$ such that
\beqtn
s_*(\dd)=\infty.
\label{butN}
\eeqtn
We argue by contradiction, and assume that for all $(\dd) \in D_{s_{0}}$, $s_*(\dd)<+\infty$. 
As we have just stated, one of the symbols '$\leq$' in the definition of $V_A(s)$ and $\tilde{V}_{\tilde{A}}(s)$ should be replaced by '$=$' symbols when $s=s_*$.\\ 
In fact, this is possible only with the '$\leq$' symbols concerning the components $q_0$, $q_1$, $\tilde{q}_0$ or $\tilde{q}_1$, as one sees in the following:

\begin{lemma}\label{transversality}{\bf (Reduction to a finite dimensional problem)}
There exists $A_3>0$ such that for each $A\geq A_3$ and $\tilde{A}\geq A_3$ there exists $s_{0,3}(A,\tilde{A})\geq s_{0,2}(A,\tilde{A})$ such that if $s_0\geq s_{0,3}$, then $(q_0(s_*),q_1(s_*),\tilde{q}_0(s_*),\tilde{q}_1(s_*))\in \pa\left( \left[-\frac{A}{s_{*}^{2}},\frac{A}{s_{*}^{2}}\right]^2\times  \left[-\frac{\tilde{A}}{s_{*}^{\alpha}},\frac{\tilde{A}}{s_{*}^{\alpha}}\right]^{2}\right)$.
\end{lemma}
\textit{Proof}: This is a direct consequence of the dynamics of equation (\ref{qqtilde}), as we will show in Subsection \ref{reduction} below.\\
Just to give a flavor of the argument, we invite the reader to look at Proposition \ref{prop36} below, where we project system (\ref{qqtilde}) on the different components of $q$ and $\tilde{q}$ introduced in (\ref{decompq}). There, one can see that the components $q_2$, $q_-$ and $q_e$ (respectively $\tilde{q}_2$, $\tilde{q}_-$ and $\tilde{q}_e$) correspond to decreasing directions of the flow and since they are ``small`` at $s=s_0$
 (see Lemma \ref{prop37} below), 
they remain small up for $s\in [s_0,s_*]$, and can not touch their bounds. Thus, only $q_0$, $q_1$, $\tilde{q}_0$ or $\tilde{q}_1$ may touch their boundary at $s=s_*$.\\
For more details on the arguments, see Subsection \ref{reduction} below. This ends the proof of  Lemma \ref{transversality}.$\blacksquare$

\bigskip

From Lemma \ref{transversality}, we may define the rescaled flow $\Phi$ at $s=s_*$ for the four expanding directions, namely $q_0$, $q_1$, $\tilde{q}_0$ and $\tilde{q}_1$, as follows:
\beqtn
\begin{array}{lccl}
\Phi : & D_{s_0}&\to& \pa([-1,1]^4)\\
             &(d_0,d_1,\tilde{d}_0,\tilde{d}_1)&\to& \displaystyle \left(\frac{s_{*}^{2}q_0}{A},\frac{s_{*}^{2}q_1}{A},\frac{s_{*}^{\alpha}\tilde{q}_0}{\tilde{A}},\frac{s_{*}^{\alpha}\tilde{q}_1}{\tilde{A}}\right)_{d_0,d_1,\tilde{d}_0,\tilde{d}_1}(s_*) .    
\end{array}
\label{defphi}
\eeqtn
In particular, 
\beqtn
\mbox{either }\omega q_m(s_*)=\frac{A}{s_{*}^{2}}\mbox{ or }\tilde{\omega} q_{\tilde{m}}(s_*)=\frac{\tilde{A}}{s_{*}^{\alpha}}, 
\label{qmqtildem}
\eeqtn
for some $m$, $\tilde{m}\in \{0,1\}$ and  $\omega$, $\tilde{\omega}\in \{-1,1\}$, both depending on
$(\dd)$. 
In the following lemma, we show that $q_m$ (or $\tilde{q}_m$) actually crosses its boundary at $s=s_*$, resulting in the continuity of $s_*$ and $\Phi$. More precisely, we have the following:
\begin{lemma}[Transverse crossing] \label{transcross}
 There exists $\Ag>0$ such that for all $A\ge \Ag$ and $\tilde A \ge \Ag$, there exists $\sg(A, \tilde A)\ge \sf(A, \tilde A)$ such that if $s_0 \ge \sg$ and (\ref{qmqtildem}) holds for some $s_*\ge \sg$, then 
 \beqtn
\mbox{either }\omega\frac{d q}{ds}(s_*)>0\mbox{ or }\tilde{\omega}\frac{d \tilde{q}_{\tilde{m}}}{ds}(s_*)>0.
\label{qmqtildemin}
\eeqtn
\end{lemma}
Clearly, from the transverse crossing, we see that
 \[
(\dd)\mapsto s_*(\dd) 
\mbox{ is continuous,}
\]
 hence by definition (\ref{defphi}), $\Phi$ is continuous.
In order to find a contradiction and conclude, we crucially use the particular form we choose for initial data in (\ref{initialq}). More precisely, we have the following:
\begin{lemma}[Degree $1$ on the boundary]\label{lem3}
There exists $\Ah>0$ such that for each $A\geq \Ah$ and $\tilde{A}\geq \Ah$, there exists $\sh(A,\tilde{A})\geq \sg(A,\tilde{A})$ such that if $s_0\geq \sg$, then the mapping $(d_0,d_1,\tilde{d}_0,\tilde{d}_1)\to (q_0(s_0),q_1(s_0),\tilde{q}_0(s_0),\tilde{q}_1(s_0))$ maps $\pa D_{s_0}$ into $\pa\left ( [-\frac{A}{s_{0}^{2}},\frac{A}{s_{0}^{2}}]^2\times [-\frac{\tilde{A}}{s_{0}^{\alpha}},\frac{\tilde{A}}{s_{0}^{\alpha}}]^2\right)$, and has degree one on the boundary.
\end{lemma}
 Indeed, from this lemma and the transverse crossing property of Lemma \ref{transcross}, we see that if $(\dd)\in\pa D_{s_0}$, then
 $s_*(\dd)=s_0$,  
$\Phi(s_*(\dd),\dd)=\left(\frac{s^{2}_{*}q_0}{A},\frac{s^{2}_{*}q_1}{A},\frac{s^{\alpha}_{*}\tilde{q}_0}{A},\frac{s^{\alpha}_{*}\tilde{q}_1}{A}\right)(s_0)$ and $\Phi$ defined in (\ref{defphi}) is a continuous function from the cuboid $D_{s_0}\subset \R^4$ to $\pa [-1,1]^4$, whose restriction to $\pa D_{s_0}$ is of degree 1. This is a contradiction. Thus, Proposition \ref{Newprop} is proved, and from identity (\ref{petit}), we have constructed a solution $(q,\tilde{q})$ to system (\ref{qqtilde}), such that
\[\|q(s)\|_{L^\infty}+\|\tilde{q}(s)\|_{L^\infty}\to 0\mbox{ as }s\to\infty.\]
In the following subsections, we give the proofs of the technical steps of the current subsection (namely Lemmas \ref{initialisationN}, \ref{transversality}, \ref{transcross} and \ref{lem3}), referring to earlier work when the proof is the same.

 \subsection{Preparation of initial data}
 In this subsection, we study initial data given by (\ref{initialq}). More precisely, we state a lemma which directly implies Lemmas \ref{initialisationN} and \ref{lem3}. It also shows the (relative) smallness of the components $q_2$, $q_-$, $q_e$, $\tilde{q}_2$, $\tilde{q}_-$ and $\tilde{q}_e$, an information which will be useful for the next subsection, dedicated to the dynamics of equation (\ref{qqtilde}), crucial for the proofs of the reduction to a finite dimensional problem (Lemma \ref{transversality}) and the transverse crossing property (Lemma \ref{transcross}). More precisely, we claim the following:
\begin{lemma}\label{prop37}{\bf (Decomposition of initial data in different components)}
For each $|\tilde{B}|\leq 1$, $A\geq 1$, $\tilde{A}\geq 1$, $0<\eta<\frac{1}{10}$ and $2<\alpha<2+\eta$, there exists 
$\si(A,\tilde{A},\tilde{B})\geq e$ 
such that for all $s_0\geq \si$:\\
 (i) there exists a cuboid
\beqtn
D_{s_0} \subset [-2,2]^4,
\eeqtn
 such that the mapping $(d_0,d_1,\tilde{d}_0,\tilde{d}_1)\to (q_0(s_0),q_1(s_0),\tilde{q}_0(s_0),\tilde{q}_1(s_0))$ is linear and one to one from $D_{s_0}$ onto $[-\frac{A}{s_{0}^{2-\eta}},\frac{A}{s_{0}^{2-\eta}}]^2\times [-\frac{\tilde{A}}{s_{0}^{\alpha}},\frac{\tilde{A}}{s_{0}^{\alpha}}]^2$ and maps the boundary $\pa D_{s_0}$ into the boundary $\pa\left ( [-\frac{A}{s_{0}^{2}},\frac{A}{s_{0}^{2}}]^2\times [-\frac{\tilde{A}}{s_{0}^{\alpha}},\frac{\tilde{A}}{s_{0}^{\alpha}}]^2\right)$. Moreover, it is of degree one on the boundary.\\
 (ii) For all $(d_0,d_1,\tilde{d}_0,\tilde{d}_1)\in D_{s_0}$, we have
 \beqtn
 \begin{array}{l}
 |q_2(s_0)|\leq C A e^{-\gamma s_0}\mbox{, for some }\gamma>0,\;\; |q_-(y,s_0)|\leq \frac{c}{s_{0}^{2}}(1+|y|^3)\mbox{ and }q_e(y,s_0)= 0,\\
 |d_0|+|d_1|\leq 1,
 \end{array}
 \eeqtn
and
\beqtn
 \begin{array}{l}
 |\tilde{q}_2(s_0)-\frac{\tilde{B}}{s_{0}^{2}}|\leq C \tilde{A} e^{-\gamma s_0} \mbox{, for some }\gamma>0,\;\; |\tilde{q}_-(y,s_0)|\leq \frac{c}{s_{0}^{\alpha}}(1+|y|^3)\mbox{ and }\tilde{q}_e(y,s_0)= 0,\\
 |\tilde{d}_0|+|\tilde{d}_1|\leq 1.
 \end{array}
 \eeqtn
\end{lemma}
{\textit{Proof}}: Since we have almost the same definition of the set $V_A$, and almost the same expression of initial data as in \cite{MZDuke97}, we refer the reader to Lemma 3.5 page 156 and Lemma 3.9 page 160 from \cite{MZDuke97}. $\blacksquare$

\subsection{Reduction to a finite-dimensional problem}
\label{reduction}
 This subsection is dedicated to the proof of Lemmas \ref{transversality} and \ref{transcross}. They both follow from the understanding of the flow of equation (\ref{qqtilde}) in the set $V_A(s)\times\tilde{V}_{\tilde{A}}(s)$. Accordingly, this crucially relies on the projection of equation (\ref{qqtilde}) with respect to the decomposition given in (\ref{decompq}). More precisely, we claim the following:
\begin{prop}[Dynamics of equation (\ref{qqtilde})]\label{prop36}
There exists $\Aj\geq 1 $ such that for all $A\geq \Aj$, $\tilde{A}\geq \Aj$,  $0<\eta<\frac{1}{10}$, $2<\alpha<2+\eta$ and $\theta\geq 0$, there exists $\sj(A,\tilde{A},\theta)$ 
such that the following holds for all $s_0\geq \sj$:\\
 Assume that for some $\tau\geq s_0$ and for all $s\in [\tau,\tau+\theta]$,
 \[(q(s),\tilde{q}(s))\in V_A(s)\times\tilde{V}_A(s).\]
Then, the following holds for all $s\in [\tau,\tau+\theta]$:\\
(i)(Differential inequalities satisfied by the expanding and null modes) For $m=0$ and $1$, we have:
\[\left |q'_m(s)-(1-\frac{m}{2})q_m(s)\right |\leq \frac{C}{s^2},\]
\[\left |\tilde{q}'_m(s)-(1-\frac{m}{2})\tilde{q}_m(s)\right |\leq \frac{C\tilde{A}^2}{s^{3-\eta}},\]
\[\left |\tilde{q}^{'}_{2}(s)+\frac{2}{s}\tilde{q}_2(s)\right |\leq \frac{C\tilde{A}}{s^{\alpha+1}}.\]
(ii)(Control of the null and negative modes) Moreover, we have:
\begin{align*}
|q_2(s)|& \leq  \frac{\tau^2}{s^2}|q_2(\tau)|+\frac{C A  (s-\tau)}{s^3},\\
|\tilde{q}_{2}(s)|& \leq  \frac{\tau^2}{s^2}|\tilde{q}_2(\tau)|+\frac{C\tilde{A}(s-\tau)}{s^{\alpha+1}},\\
\left\| \frac{q_-(s)}{1+|y|^3}\right\|_{L^{\infty}}&\leq C e^{-\frac{(s-\tau)}{2}} \left\| \frac{q_-(\tau)}{1+|y|^3}\right\|_{L^{\infty}}
+C\frac{e^{-(s-\tau)^2}\|q_e(\tau)\|_{L^\infty}}{s^{3/2}}+\frac{C(1+s-\tau)}{s^2},\\
\left\| \frac{\tilde{q}_-(s)}{1+|y|^3}\right\|_{L^{\infty}}&\leq C e^{-\frac{(s-\tau)}{2}} \left\| \frac{\tilde{q}_-(\tau)}{1+|y|^3}\right\|_{L^{\infty}}
+C\frac{e^{-(s-\tau)^2}\|\tilde{q}_e(\tau)\|_{L^\infty}}{s^{3/2}}+\frac{C(1+s-\tau)}{s^\alpha},
\\
\|q_e(s)\|_{L^\infty}&\leq C e^{-\frac{(s-\tau)}{2}}\|q_e(\tau)\|_{L^\infty}+C\frac{e^{s-\tau}}{s^{3/2}} \left\| \frac{q_-(\tau)}{1+|y|^3}\right\|_{L^\infty}+\frac{C(1+s-\tau)}{s^{1/2}},\\
\|\tilde{q}_e(s)\|_{L^\infty}&\leq C e^{-\frac{(s-\tau)}{2}}\|\tilde{q}_e(\tau)\|_{L^\infty}+C \frac{e^{s-\tau}}{ s^{3/2}} \left\| \frac{\tilde{q}_-(\tau)}{1+|y|^3}\right\|_{L^\infty}+\frac{C(1+s-\tau)}{s^{\alpha-3/2}}.
\end{align*}
\end{prop}
Let us first insist on the fact that the derivation of Lemmas \ref{transversality} and \ref{transcross} follows from Proposition \ref{prop36}, exactly as in the real case treated in \cite{MZDuke97} (see pages 163 to 166 and 158 to 159 in \cite{MZDuke97} ). For that reason, we only focus in the following on the proof of Proposition \ref{prop36}.

\medskip

\textit{Proof of Proposition \ref{prop36}}:

\medskip

\noindent The proof of Proposition \ref{prop36} consists in the projection of the two equations of system (\ref{qqtilde}) on the different components of $q$ and $\tilde{q}$ defined in (\ref{decompq}).\\
When $\tilde{q}\equiv 0$, the proof is already available from Lemma 3.13 pages 167 and Lemma 3.8 page 158 from \cite{MZDuke97}.\\
When $\tilde{q}\not\equiv 0$, since the equation satisfied by $\tilde{q}$ in (\ref{qqtilde}) shares the same linear part as the equation in $q$, the proof is similar to the argument in \cite{MZDuke97}. For that reason, we only give the ideas here, and kindly ask the interested reader to look at Lemma 3.13 page 167 and Lemma 3.8 page 158 in \cite{MZDuke97} for the technical details.

\bigskip

(i) Multiplying the two equations in (\ref{qqtilde}) by $\chi(y,s)k_m(y)\rho(y)$, for $m=0,1,2$ and integrating in $y\in\R$, we proceed as in pages 158-159 from \cite{MZDuke97} and we get the differential inequalities given in (i) with no difficulties.

\bigskip

(ii) For convenience, we separate the contribution of $q$ and $\tilde{q}$ in the quadratic term $b(y,s)$ defined in (\ref{defN}) by writing $b(y,s)=B(y,s)-N(y,s)$ with 
\beqtn
B(y,s)=q^2\mbox{, and }N(y,s)=\tilde{q}^2.\label{defNN} 
\eeqtn
Let us first recall equations of $(q,\tilde{q})$ in their Duhamel formulation,
\beqtn
\begin{array}{lll}
q(s)&=&K(s,\tau)q(\tau)+\int_{\tau}^{s}d\sigma K(s,\sigma)B(q(\sigma))+\int_{\tau}^{s}d\sigma K(s,\sigma)R(\sigma)-\int_{\tau}^{s}d\sigma K(s,\sigma)N(\sigma),\\
\tilde{q}(s)&=&K(s,\tau)\tilde{q}(\tau)+\int_{\tau}^{s}d\sigma K(s,\sigma)\tilde{b}(\sigma),
\end{array}
\eeqtn
where $K$ is the fundamental solution of the operator $\Lg+V$.
We write $q=\alpha+\beta+\gamma+\delta$ and $\tilde{q}=\tilde{\alpha}+\tilde{\beta}$, where
\beqtn
\begin{array}{l}
\alpha(s)=K(s,\tau)q(\tau),\;\;\beta(s)=\int_{\tau}^{s}d\sigma K(s,\sigma)B(q(\sigma)),\\
\gamma(s)=\int_{\tau}^{s}d\sigma K(s,\sigma)R(\sigma),\;\;\delta(s)=-\int_{\tau}^{s}d\sigma K(s,\sigma)N(\sigma) \label{defdelta}.
\end{array}
\eeqtn
\beqtn
\tilde{\alpha}(s)=K(s,\tau)\tilde{q}(\tau),\;\;\tilde{\beta}(s)=\int_{\tau}^{s}d\sigma K(s,\sigma)\tilde{b}(\sigma).\label{defbeta}\\
\eeqtn
We assume that $(q(s),\tilde{q}(s))\in V_A(s)\times\tilde{V}_A(s)$ for each $s\in [\tau,\tau+\theta]$. Clearly, proceeding as the derivation of Lemma 3.13 page 167 in \cite{MZDuke97}, (ii) of Proposition \ref{prop36} follows from  the following:
\begin{lemma}[Projection of the Duhamel formulation] \label{lemduh}There exists $\Ak\geq 1$ such that for all $A\geq \Ak$, $\tilde{A}\geq \Ak$ and $\theta>0$ there exists $\sk(A,\tilde{A},\theta)\geq \sj(A)$, such that for all $s_0\geq \sk(A,\tilde{A},\theta)$, if we assume that for some $\tau\geq s_0$ and for all $s\in[\tau,\tau+\theta]$, $q(s)\in V_A(s)$ and $\tilde{q}(s)\in\tilde{V}_{\tilde{A}}(s)$, then\\
(i) (Linear terms)
\beqtn
\begin{array}{lll}
|\alpha_2(s)|&\leq& \frac{\tau^2}{s^2}|q_2(\tau)|+\frac{C A  (s-\tau)}{s^3},\\
\left\| \frac{\alpha_-(s)}{1+|y|^3}\right\|_{L^{\infty}}&\leq &C e^{-\frac{(s-\tau)}{2}} \left\| \frac{q_-(\tau)}{1+|y|^3}\right\|_{L^{\infty}}
+C\frac{e^{-(s-\tau)^2}\|q_e(\tau)\|_{L^\infty}}{s^{3/2}}+\frac{C}{s^2},\\
\|\alpha_e(s)\|_{L^\infty}&\leq &Ce^{-\frac{(s-\tau)}{2}}\|q_e(\tau)\|_{L^\infty}+Ce^{s-\tau} s^{3/2} \left\| \frac{q_-(\tau)}{1+|y|^3}\right\|_{L^\infty}+\frac{C}{\sqrt{s}},
\end{array}
\eeqtn
and
\beqtn
\begin{array}{lll}
|\tilde{\alpha}_2(s)|&\leq &\frac{\tau^2}{s^2}|\tilde{q}_2(\tau)|+\frac{C A  (s-\tau)}{s^{\alpha+1}},\\
\left\| \frac{\tilde{\alpha}_-(s)}{1+|y|^3}\right\|_{L^{\infty}}&\leq& C e^{-\frac{(s-\tau)}{2}} \left\| \frac{\tilde{q}_-(\tau)}{1+|y|^3}\right\|_{L^{\infty}}
+C\frac{e^{-(s-\tau)^2}\|\tilde{q}_e(\tau)\|_{L^\infty}}{s^{3/2}}+\frac{C}{s^\alpha},\\
\|\tilde{\alpha}_e(s)\|_{L^\infty}&\leq & Ce^{-\frac{(s-\tau)}{2}}\|\tilde{q}_e(\tau)\|_{L^\infty}+Ce^{s-\tau} s^{3/2} \left\| \frac{\tilde{q}_-(\tau)}{1+|y|^3}\right\|_{L^\infty}+\frac{C}{s^{\alpha-3/2}}.
\end{array}
\eeqtn
(ii) (Nonlinear terms)
\begin{align*}
|\beta_2(s)|&\leq \frac{(s-\tau)}{s^3},&|\beta_-(y,s)|&\leq \frac{(s-\tau)}{s^2}(1+|y|^3),&\|\beta_e(s)\|_{L^\infty}&\leq \frac{(s-\tau)}{\sqrt s},\\
|\delta_2(s)|&\leq C\frac{(s-\tau)}{s^3}, &|\delta_-(y,s)|&\leq C\frac{(s-\tau)}{s^2}(1+|y|^3),&\|\delta_e(s)\|_{L^\infty}&\leq C\frac{(s-\tau)}{\sqrt s},\\
|\tilde{\beta}_2(s)|&\leq \frac{(s-\tau)}{s^{\alpha+1}},&|\tilde{\beta}_-(y,s)|&\leq \frac{(s-\tau)}{s^\alpha}(1+|y|^3),&\|\tilde{\beta}_e(s)\|_{L^\infty}&\leq \frac{(s-\tau)}{s^{\alpha-3/2}}.
\end{align*}
(iii) (Source term)
\[|\gamma_2(s)|\leq C(s-\tau)s^{-3},\;|\gamma_-(y,s)|\leq C(s-\tau)(1+|y|^3)s^{-2},\;\|\gamma_e(s)\|_{L^\infty}\leq (s-\tau)s^{-1/2}.\]
\end{lemma}
\textit{Proof:} We consider, $A\geq 1$, $\tilde{A}\geq 1$, $\theta>0$, 
and $s_0\ge \theta$. 
The terms $\alpha$, $\beta$ and $\gamma$ are already present in the case of the real-valued semilinear heat equation, so we refer to Lemma 3.13 page 167 in \cite{MZDuke97} for the estimates involving them. As for $\tilde{\alpha}$, since the definition of $\tilde{V}_{\tilde{A}}(s)$ is different from the definition of $V_A(s)$, the reader will have absolutely no difficulty to adapt Lemma 3.13 of \cite{MZDuke97} to the new situation. Thus, we only focus on the new terms $\delta(y,s)$ and $\tilde{\beta}(y,s)$. Note that since $s_0\geq \theta$, if we take $\tau\geq s_0$, then $\tau+\theta\leq 2 \tau$ and if $\tau \leq \sigma\leq s\leq \tau +\theta$, then
\[\frac{1}{2\tau}\leq \frac{1}{s}\leq \frac{1}{\sigma}\leq \frac{1}{\tau}.\]
Let us first derive the following bounds when $(q(s),\tilde{q}(s))\in V_A(s)\times\tilde{V}_{\tilde{A}}(s)$:
\begin{prop} [Bounds for $(q(s),\tilde{q}(s))\in V_A(s)\times\tilde{V}_{\tilde{A}}(s)$]
 \label{VA}
For all $s\geq e$, we consider $r\in V_A(s)$ and $\tilde{r}\in \tilde{V}_{\tilde{A}}(s)$, where the shrinking sets $V_A(s)$ and $ \tilde{V}_{\tilde{A}}(s)$ are given in Definition \ref{defVA}.  Then, we have:
 \[\begin{array}{l}
 (i)\; for \;all\; y \in \R,\;\; |r(y,s)|\leq C A^2 \frac{\log s}{s^2}(1+|y|^3),\\
(ii)\;  \|r(s)\|_{L^\infty}\leq C\frac{A^2}{\sqrt{s}},\\
(iii)\; for\; all\; y \in \R,\;\; |\tilde{r}(y,s)|\leq C \frac{\tilde{A}^2}{s^{2-\eta}}(1+|y|^3),\;\;|\tilde{r}_b(y,s)|\leq \frac{C\tilde{A}}{s^{\alpha-3/2}},\\
(iv)\;  \|\tilde{r}(s)\|_{L^\infty}\leq C\frac{\tilde{A}^2}{s^{\alpha-3/2}}.
\end{array}
\]
 \end{prop}
 \textit{Proof}: The proof is omitted since it is the same as the corresponding part in \cite{MZDuke97}. See Proposition 3.7 page 157 in   \cite{MZDuke97} for details. $\blacksquare$

\bigskip

Then, we recall from Bricmont and Kupiainen \cite{BKN94} the following estimates on $K(s,\sigma)$, the semigroup generated by $\Lg+V$: 
\begin{lemma}\label{lemK}(Properties of $K(s,\sigma)$):\\
(i) For all $s\geq \sigma>1$ and $y$, $x\in\R$, we have
\[|K(s,\sigma,y,x)|\leq C e^{(s-\sigma)\Lg}(y,x),\]
where $e^{\psi \Lg}$ is given by
\[e^{\psi \Lg}(y,x)=\frac{e^\psi}{\sqrt{4\pi(1-e^{-\psi})}} \exp\left[-\frac{(ye^{-\psi/2}-x)^2}{4(1-e^{-\psi})}\right].\]
(ii)We have for all $s\geq \tau\geq 1$, with $s\leq 2\tau $,
\beqtn
\displaystyle\left| \int K(s,\tau,y,x)(1+|x|^m)dx \right |\leq C\int e^{(s-\tau)\Lg}(y,x) (1+|x|^m)dx \leq e^{s-\tau} (1+|y|^m).
\eeqtn
\end{lemma}
{\it Proof}:\\
(i) See page 181 in \cite{MZDuke97}\\
(ii) See Corollary 3.14 page 168 in \cite{MZDuke97}. $\blacksquare$
 
 \medskip
 
 \noindent{\it Estimates on $\delta$ defined in (\ref{defdelta})}:\\
 Consider $s\in[\tau,\tau+\theta]$. Since $\tilde{q}(s)\in\tilde{V}_{A}(s)$ by assumption, using (iii) and (iv) of Lemma \ref{VA}, we see that 
 \beqtn
 \forall y\in\R,\;|\tilde{q}(y,s)|\leq \min\left( \frac{C\tilde{A}^2}{s^{2-\eta}}(1+|y|^3),\frac{C\tilde{A}^2}{s^{\alpha-3/2}}\right),
 \label{bqt}
 \eeqtn
hence by definition (\ref{defNN}) of $N$, we obtain
\beqtn
\forall y \in\R,\; |N(y,s)|\leq C \tilde{A}^4\min\left( \frac{(1+|y|^3)}{s^{\alpha+\frac 1 2-\eta}}, \frac{1}{s^{2\alpha-3}},\frac{1+|y|^6}{s^{4-2\eta}}\right).
\label{bN}
\eeqtn
Using Lemma \ref{lemK} and the  definition (\ref{defdelta}) of $\delta$, we write
\beqtn
\begin{array}{lll}
|\delta(y,s)|  &\leq&\displaystyle\int_{\tau}^{s} d\sigma\int_{\R}\left|K(s,\sigma,y,x)N(x,\sigma) \right|dx\\
&\leq&\displaystyle \int_{\tau}^{s} d\sigma\int_{\R} e^{(s-\sigma)\Lg}(y,x)\frac{C\tilde{A}^4(1+|x|^3)}{s^{\alpha+1/2-\eta}} dx\\
&\leq &\displaystyle\frac{C\tilde{A}^4 (s-\tau)}{s^{\alpha+1/2-\eta}} e^{s-\tau}(1+|y|^3)\leq \frac{(s-\tau)}{s^2}(1+|y|^3),\\
\label{esdelta}
\end{array}
\eeqtn
for $s_0$ large enough, since $\eta<1/2$.\\
Using the following bounds in (\ref{bN}) and proceeding similarly, we see that
\[\forall y\in\R,\;|\delta(y,s)|\leq (s-\tau)\min\left(\frac{1+|y|^3}{s^2},\frac{1}{\sqrt{s}},\frac{1+|y|^6}{s^3}\right),\] 
since $\alpha> 2$, $\eta<1/2$, and provided that $s_0$ is large enough.\\
By definition of $q_m$, $q_-$ and $q_e$ for $m\leq 2$, we write

\beqtn
\begin{array}{lll}
|\delta_m(s)|&\leq&\left|\int_{\R}\chi(y,s)\delta(y,s)k_m(y)\rho(y)dy \right|\leq C\int_{\R}|\delta(y,s)|(1+|y|^2)\rho(y)dy\leq\frac{ C(s-\tau)}{s^{3}},\\
|\delta_-(y,s)|&=&\left|\chi(y,s) \delta(y,s)-\sum_{i=0}^{2}\delta_i(s)k_i(y)\right|\leq(s-\tau)(1+|y|^3)\frac{C}{s^{2}}.\\
|\delta_e(y,s)|&=&\left|(1-\chi(y,s))\delta (y,s)\right\|\leq(s-\tau)\frac{C}{\sqrt{s}}.
\label{details}
\end{array}
\eeqtn
\noindent{\it Estimates on $\tilde{\beta}$ defined in (\ref{defbeta}):}\\
Consider $s\in[\tau,\tau+\theta]$. Since $q(s)\in V_A(s)$ by assumption, using (i) and (ii) of Lemma \ref{VA}, we see that
\[\forall y\in\R,\; |q(y,s)|\leq CA^2 \min\left( \frac{\log s}{s^2}(1+|y|^3),\frac{1}{\sqrt{s}}\right).\]
Using (\ref{bqt}) and the definition (\ref{defN}) of $\tilde{b}$, we see that
\[\forall y\in\R,\; |\tilde{b}(y,s)|\leq CA^2\tilde{A}^2 \min\left( \frac{\log s}{s^{\alpha+1/2}}(1+|y|^3),\frac{1}{s^{\alpha-1}}, \frac{1+|y|^6}{s^{4-\eta}}\right).\]
Using the definition (\ref{defbeta}) of $\tilde{\beta}$ and arguing as for estimate (\ref{esdelta}), we see that
\[\forall y\in\R,\; |\tilde{\beta}(y,s)|\leq (s-\tau)\min\left( \frac{1+|y|^3}{s^\alpha},\frac{1}{s^{\alpha-3/2}},\frac{1+|y|^6}{s^{\alpha+1}}\right),\]
provided that $s_0$ is large enough, since $\eta<\frac 1 2$ and $\alpha<2+\eta<3-\eta$. Arguing as for (\ref{details}), we get the desired estimates. This concludes the proof of Lemma \ref{lemduh}. $\blacksquare$

\medskip

Since item (ii) of Proposition \ref{prop36} follows from Lemma \ref{lemduh}, exactly as for Lemma 3.13 page 167 in \cite{MZDuke97}, this also ends the proof of Proposition \ref{prop36} . $\blacksquare$

\medskip

Since Lemmas \ref{transversality} and \ref{transcross} follow from Proposition \ref{prop36} exactly as in \cite{MZDuke97} (see pages 163 to 166 and 158 to 159 in that paper), this is also the conclusion of the proof of Lemmas \ref{transversality} and \ref{transcross}.
Recalling that we have already justified that Lemmas \ref{initialisationN} and \ref{lem3} hold (see Lemma \ref{prop37} above), and given that Proposition \ref{Newprop} is the consequence of Lemmas \ref{initialisationN}, \ref{transversality}, \ref{transcross} and \ref{lem3}, this is also the conclusion of Proposition \ref{Newprop}. $\blacksquare$

\section{Asymptotic behavior of $u(t)$}
We prove Theorem \ref{theorem1} in this section. We will first derive (ii) and (iii) from Section 3, then we will prove (i) and (iv). \\
Consider $0<|\tilde{B}|\leq 1$. Using Proposition \ref{Newprop}, Lemma \ref{prop37} and Proposition \ref{prop36}, we see that if $A=\tilde{A}=\max(1,\Ad,\Aj)$, $0<\eta< \frac 1{10}$, $2<\alpha<2+\eta$ and 
$T\leq \Tl(\tilde B)$ for some $\Tl(\tilde B)\le \min(\Td, \Ti, \Tj)$ where $T_i = -\log s_{0,i}$,  
then there exists a parameter $(d_0,d_1,\tilde{d}_0,\tilde{d}_1)$ such that if $(q(s_0),\tilde{q}(s_0))$ is given by \eqref{initialq}, where $s_0=-\log T$, then 
\beqtn
\begin{array}{ll}
\forall s\geq -\log T,\; q(s)\in V_A(s),\; \tilde{q}(s)\in \tilde{V}_{\tilde{A}}(s),&\; \left| \tilde{q}_{2}^{'}(s)+\frac 2 s \tilde{q}_2(s) \right|\leq C\frac{\tilde A}{s^{\alpha+1}}
\leq\frac{\mu_0}{s^{\alpha+1}},\\
\mbox{with }\mu_0=\frac{\alpha-2}{4}|\tilde{B} |s_{0}^{\alpha-2},&
\end{array}
\label{but}
\eeqtn
and 
\[\left|\tilde{q}_2(s_0)-\frac{\tilde{B}}{s_{0}^{2}} \right|
\leq C\tilde A e^{-\gamma s_0}
\leq \frac{|\tilde{B}|}{4 s_{0}^{2}}.\]
As announced  earlier, we use this property to derive (ii) and (iii) of Theorem \ref{theorem1}, then we will prove (i) and (iv).

\medskip

\noindent (ii) This directly follows from \eqref{but}, \eqref{petit}
and the selfsimilar transformation \eqref{chauto}.\\
(iii) From \eqref{but}, we see that 
\beqtn
\forall s\geq -\log T,\;|\left(s^2\tilde{q}_2\right)^{'}|\leq \frac{\mu_0}{s^{\alpha-1}},
\label{integq2}
\eeqtn
which means that $s^2\tilde{q}_2(s)$ has some limit $l$ as $s\to\infty$.\\
Integrating this inequality between $s$ and $+\infty$, we obtain
\beqtn
|s^2\tilde{q}_2(s)-l|\leq \frac{\mu_0}{(2-\alpha)s^{\alpha-2}}.
\label{limite}
\eeqtn
Putting $s=s_0$ in this identity, then using \eqref{but}, we see that 
\[|s_{0}^{2}\tilde{q}_2(s_0)-l|\leq \frac{|\tilde{B}|}{4}\mbox{ and }|s_{0}^{2}\tilde{q}_2(s_0)-\tilde{B}|\leq \frac{|\tilde{B}|}{4},\]
Thus, it follows that
\[|l-\tilde{B}|\leq \frac{|\tilde{B}|}{2}\mbox{, hence }|l|\geq \frac{|\tilde{B}|}{2}>0\mbox{ and }l\not\equiv 0.\]
We then write from the decomposition \eqref{decompq} that for all $s\geq -\log T$, $R>0$ and $|y|\leq R$, $\tilde{q}_e(y,s)=0$, hence,
\[ \tilde{q}(y,s)-\frac{l}{s^2}h_2(y)=\sum_{i=0}^{1}\tilde{q}_i(s)h_i(y)+(\tilde{q}_2(s)-\frac{l}{s^2})h_2(y)+\tilde{q}_-(y,s).\]
Using the fact that for all $s\geq -\log T$, $\tilde{q}(s)\in\tilde{V}_{\tilde{A}}(s)$ (see \eqref{but} above), Definition \ref{defVA} for $\tilde{V}_{\tilde{A}}(s)$, together with \eqref{limite}, we see that for all $s\geq -\log T,\;R>0\mbox{ and }|y|\leq R$
\[
\left|\tilde{q}(y,s)-\frac{l}{s^2} h_2(y)\right|\le \frac{C(\tilde A, R)}{s^\alpha}.
\]
Using the definition \eqref{initialqqtilde} of $\tilde{q}$ and \eqref{chauto} of $\tilde{w}$, we get the desired conclusion.

\medskip

(i) If $x_0=0$, then we see from (\ref{profilev}) and (\ref{profilevt}) that $|v(0,t)| \sim (T-t)^{-1}$ as $t\to T$. Hence $u$ blows up at time $T$ at $x_0=0$. It remains to prove that any $a\neq 0$ is not a blow-up point. The following result from Giga and Kohn \cite{GKCPAM89} allows us to conclude:
\begin{prop}(Giga and Kohn - No blow-up under the ODE threshold) For all $C_0>0$, there is $\eta_0>0$ such that if $v(\xi,\tau)$ solves
\[\left |  v_t -\Delta v\right |\leq C_0 (1+|v|^p)\]
and satisfies
\[|v(\xi,\tau)|\leq \eta_0(T-t)^{-1}\]
for all $(\xi,\tau)\in B(a,r)\times[T-r^2,T)$ for some $a\in \R$ and $r>0$, then $v$ does not blow up at (a,T).
\label{propositionGK}
\end{prop}

\textit{Proof: } See Theorem 2.1 page 850 in \cite{GKCPAM89}. Note that the proof of Giga and Kohn is valid also when $u$ is complex valued. $\blacksquare$\\
Indeed, since we see from (\ref{profilev}) that
\[\sup_{|x-x_0|\leq |x_0|/2}(T-t)^{-1}|u(x,t)|\leq \left | \varphi\left(\frac{|x_0|/2}{\sqrt{(T-t)|\log(T-t)|}}\right)\right |+\frac{C}{\sqrt{|\log(T-t)|}}\to 0\]
as $t\to T$, $x_0$ is not a blow-up point of $u$ from Proposition \ref{propositionGK}. This concludes the proof of (i) of Theorem \ref{theorem1}.\\
(iv) Arguing as Merle did in \cite{FMCPAM92}, we derive the existence of a blow-up profile $u^*\in C^2(\R^*)$ such that $u(x,t)\to u^*(x)$ as $t\to T$, uniformly on compact sets of $\R^*$. The profile $u^*(x)$ is not defined at the origin. In the following, we would like to find its equivalent as $x\to 0$ and show that it is in 444 singular at the origin. We argue as in Masmoudi and Zaag \cite{MZ07}. Consider $K_0>0$ to be fixed large enough later. If $x_0\neq 0$ is small enough, we introduce for all $(\xi,\tau)\in \R\times [-\frac{t_0(x_0)}{T-t_0(x_0)},1)$, 

\begin{align}
V(x_0,\xi,\tau)&=(T-t_0(x_0)) v(x,t),\\
\tilde{V}(x_0,\xi,\tau)&=(T-t_0(x_0))\tilde{ v}(x,t),\\
\mbox{where,   }x&=x_0+\xi\sqrt{T-t_0(x_0)},\; t=t_0(x_0)+\tau(T-t_0(x_0)),
\label{defV}
\end{align}

and $t_0(x_0)$ is uniquely determined by
\beqtn
|x_0|=K_0\sqrt{(T-t_0(x_0))|\log(T-t_0(x_0))|}.
\label{xt0}
\eeqtn
From the invariance of problem (\ref{uv}) under dilation, $(V(x_0,\xi,\tau),\tilde{V}(x_0,\xi,\tau))$ is also a solution of (\ref{uv}) on its domain. From (\ref{defV}), (\ref{xt0}), (\ref{profilevt}) and (\ref{profilev}), we have
\[\sup_{|\xi|<2|\log(T-t_0(x_0))|^{1/4}}\left |V(x_0,\xi,0)-f(K_0) \right|\leq \frac{C}{|\log(T-t_0(x_0))|^{1/4}}\to 0\mbox{ as }x_0\to 0\]
and
\[\sup_{|\xi|<2|\log(T-t_0(x_0))|^{1/4}}\left |\tilde{V}(x_0,\xi,0)\right|\leq \frac{C}{|\log(T-t_0(x_0))|^{1/4}}\to 0\mbox{ as }x_0\to 0.\]
Using the continuity with respect to initial data for problem (\ref{uv}) associated to a space-localization in the ball $B(0,|\xi|<|\log(T-t_0(x_0))|^{1/4})$, we show as in Section 4 of \cite{ZAIHPANL98} that
\[
\begin{array}{l}
\sup_{|\xi|\leq |\log(T-t_0(x_0))|^{1/4},\;0\leq\tau<1}|V(x_0,\xi,\tau)-U_{K_0}(\tau)|\leq \epsilon(x_0)\mbox{ as }x_0\to 0 \\
\sup_{|\xi|\leq |\log(T-t_0(x_0))|^{1/4},\;0\leq\tau<1}|\tilde{V}(x_0,\xi,\tau)|\leq \epsilon(x_0)\mbox{ as }x_0\to 0, 
\end{array}
\]
where $U_{K_0}(\tau)=((1-\tau)+\frac{K_{0}^{2}}{8})^{-1}$ is the solution of the PDE (\ref{uv}) with constant initial data $\varphi(K_0)$. Making $\tau \to 1$ and using (\ref{defV}), we see that 
\[
\begin{array}{lll}
v^*(x_0)=\lim_{t\to T} v(x,t)&=&(T-t_0(x_0))^{-1}\lim_{\tau \to 1} V(x_0,0,\tau)\\
&\sim&(T-t_0(x_0))^{-1}U_{K_0}(1)
\end{array}
\]
as $x_0\to 0$. We note also that
\[|\tilde{v}^*(x_0)|\leq \epsilon(x_0) (T-t_0(x_0))^{-1}.\]
Since we have from (\ref{xt0}) 
\[\log(T-t_0(x_0)) \sim 2 \log |x_0|\mbox{ and } T-t_0(x_0)\sim \frac{|x_0|^2}{2K_{0}^{2}|\log|x_0||},\]
as $x_0\to 0$, this yields (iv) of Theorem \ref{theorem1} and concludes the proof of Theorem \ref{theorem1}. $\blacksquare$



\begin{thebibliography}{FKMZ00}

\bibitem[Bal77]{BQJMOS77}
J.~M. Ball.
\newblock Remarks on blow-up and nonexistence theorems for nonlinear evolution
  equations.
\newblock {\em Quart. J. Math. Oxford Ser. (2)}, 28(112):473--486, 1977.

\bibitem[BK94]{BKN94}
J.~Bricmont and A.~Kupiainen.
\newblock Universality in blow-up for nonlinear heat equations.
\newblock {\em Nonlinearity}, 7(2):539--575, 1994.

\bibitem[CZ13]{CZCPAM13}
R.~C{\^o}te and H.~Zaag.
\newblock Construction of a multi-soliton blow-up solution to the semilinear
  wave equation in one space dimension.
\newblock {\em Comm. Pure Appl. Math.}, 66(10):1541--1581, 2013.

\bibitem[EZ11]{EZSMJ11}
M.~A. Ebde and H.~Zaag.
\newblock Construction and stability of a blow up solution for a nonlinear heat
  equation with a gradient term.
\newblock {\em S$\vec{ \rm e}$MA J.}, 55:5--21, 2011.

\bibitem[FKMZ00]{FMZMA00}
C.~Fermanian~Kammerer, F.~Merle, and H.~Zaag.
\newblock Stability of the blow-up profile of non-linear heat equations from
  the dynamical system point of view. math.
\newblock {\em Math. Annalen}, 317(2):347--387, 2000.

\bibitem[FKZ00]{FZN00}
C.~Fermanian~Kammerer and H.~Zaag.
\newblock Boundedness up to blow-up of the difference between two solutions to
  a semilinear heat equation.
\newblock {\em Nonlinearity}, 13(4):1189--1216, 2000.

\bibitem[FM95]{FMJDE95}
S.~Filippas and F.~Merle.
\newblock Modulation theory for the blowup of vector-valued nonlinear heat
  equations.
\newblock {\em J. Differential Equations}, 116(1):119--148, 1995.

\bibitem[GK85]{GKCPAM85}
Y.~Giga and R.V. Kohn.
\newblock Asymptotically self-similar blow-up of semilinear heat equations.
\newblock {\em Comm. Pure Appl. Math.}, 38(3):297--319, 1985.

\bibitem[GK87]{GKIUMJ87}
Y.~Giga and R.~V. Kohn.
\newblock Characterizing blowup using similarity variables.
\newblock {\em Indiana Univ. Math. J.}, 36(1):1--40, 1987.

\bibitem[GK89]{GKCPAM89}
Y.~Giga and R.~V. Kohn.
\newblock Nondegeneracy of blowup for semilinear heat equations.
\newblock {\em Comm. Pure Appl. Math.}, 42(6):845--884, 1989.

\bibitem[GNSY13]{GHMY12}
J.S. Guo, H.~Ninomiya, M.~Shimojo, and Yanagida.E.
\newblock Convergence and blow-up of solutions for a complex-valued heat
  equation with a quadratic nonlinearity.
\newblock {\em Trans. Amer. Math. Soc.}, 365(5):2447--2467, 2013.

\bibitem[HV92]{HVADSSP92}
M.~A. Herrero and J.~J.~L. Vel{\'a}zquez.
\newblock Generic behaviour of one-dimensional blow up patterns.
\newblock {\em Ann. Scuola Norm. Sup. Pisa Cl. Sci. (4)}, 19(3):381--450, 1992.

\bibitem[HV93]{HVAIHP93}
M.~A. Herrero and J.~J.~L. Vel{\'a}zquez.
\newblock Blow-up behaviour of one-dimensional semilinear parabolic equations.
\newblock {\em Ann. Inst. H. Poincar\'e Anal. Non Lin\'eaire}, 10(2):131--189,
  1993.

\bibitem[HV94]{HVCRAS94}
M.A. Herrero and J.J.~L. Vel{\'a}zquez.
\newblock Explosion de solutions d'\'equations paraboliques semilin\'eaires
  supercritiques.
\newblock {\em C. R. Acad. Sci. Paris S\'er. I Math.}, 319(2):141--145, 1994.

\bibitem[Mer92]{FMCPAM92}
F.~Merle.
\newblock Solution of a nonlinear heat equation with arbitrary given blow-up
  points.
\newblock {\em Comm. Pure Appl. Math.}, 45(3):263--300, 1992.

\bibitem[Miz07]{MMA07}
N.~Mizoguchi.
\newblock Rate of type {II} blowup for a semilinear heat equation.
\newblock {\em Math. Ann.}, 339(4):839--877, 2007.

\bibitem[MM04]{MMCPAM04}
H.~Matano and F.~Merle.
\newblock On nonexistence of type {II} blowup for a supercritical nonlinear
  heat equation.
\newblock {\em Comm. Pure Appl. Math.}, 57(11):1494--1541, 2004.

\bibitem[MM09]{MMJFA08}
H.~Matano and F.. Merle.
\newblock Classification of type {I} and f{II} behaviors for a supercritical
  nonlinear heat equation.
\newblock {\em J. Funct. Anal.}, 256(4):992--1064, 2009.

\bibitem[MRR11]{MRRCRAS11}
F.~Merle, P.~Rapha{\"e}l, and I.~Rodnianski.
\newblock Blow up dynamics for smooth equivariant solutions to the energy
  critical {S}chr\"odinger map.
\newblock {\em C. R. Math. Acad. Sci. Paris}, 349(5-6):279--283, 2011.

\bibitem[MZ97]{MZDuke97}
F.~Merle and H.~Zaag.
\newblock Stability of the blow-up profile for equations of the type
  $u_t={\Delta} u+\vert u\vert ^{p-1}u$.
\newblock {\em Duke Math. J.}, 86(1):143--195, 1997.

\bibitem[MZ98]{MZCPAM98}
F.~Merle and H.~Zaag.
\newblock Optimal estimates for blowup rate and behavior for nonlinear heat
  equations.
\newblock {\em Comm. Pure Appl. Math.}, 51(2):139--196, 1998.

\bibitem[MZ00]{MZMA00}
F.~Merle and H.~Zaag.
\newblock A {L}iouville theorem for vector-valued nonlinear heat equations and
  applications.
\newblock {\em Math. Ann.}, 316(1):103--137, 2000.

\bibitem[MZ08]{MZ07}
N.~Masmoudi and H.~Zaag.
\newblock Blow-up profile for the complex {G}inzburg-{L}andau equation.
\newblock {\em J. Funct. Anal.}, 225:1613--1666, 2008.

\bibitem[OSW08]{OSWN08}
H.~Okamoto, T.~Sakajo, and M.~Wunsch.
\newblock On a generalization of the constantin-lax-majda equation.
\newblock {\em Nonlinearity}, 21(10):2447--2461, 2008.

\bibitem[QS07]{QSBV07}
P.~Quittner and P.~Souplet.
\newblock {\em Superlinear parabolic problems}.
\newblock Birkh\"auser Advanced Texts: Basler Lehrb\"ucher. [Birkh\"auser
  Advanced Texts: Basel Textbooks]. Birkh\"auser Verlag, Basel, 2007.
\newblock Blow-up, global existence and steady states.

\bibitem[RR12]{RRIHES12}
P.~Rapha{\"e}l and I.~Rodnianski.
\newblock Stable blow up dynamics for the critical co-rotational wave maps and
  equivariant {Y}ang-{M}ills problems.
\newblock {\em Publ. Math. Inst. Hautes \'Etudes Sci.}, pages 1--122, 2012.

\bibitem[RS13]{RSCPAM13}
P.~Rapha{\"e}l and R.~Schweyer.
\newblock Stable blowup dynamics for the 1-corotational energy critical
  harmonic heat flow.
\newblock {\em Comm. Pure Appl. Math.}, 66(3):414--480, 2013.

\bibitem[Sak03a]{SJMSUT03}
T.~Sakajo.
\newblock Blow-up solutions for the constantin-lax-majda equation with a
  generalized viscosity term.
\newblock {\em J. Math. Sci. Univ. Tokyo.}, 10(1):187--207, 2003.

\bibitem[Sak03b]{SN03}
T.~Sakajo.
\newblock On global solutions for the constantin-lax-majda equation with a
  generalized viscosity term.
\newblock {\em Nonlinearity.}, 16:1319--1328, 2003.

\bibitem[Sch12]{SJFA12}
R.~Schweyer.
\newblock Type {II} blow-up for the four dimensional energy critical semi
  linear heat equation.
\newblock {\em J. Funct. Anal.}, 263(12):3922--3983, 2012.

\bibitem[Vel92]{VCPDE92}
J.~J.~L. Vel{\'a}zquez.
\newblock Higher-dimensional blow up for semilinear parabolic equations.
\newblock {\em Comm. Partial Differential Equations}, 17(9-10):1567--1596,
  1992.

\bibitem[Vel93a]{VTAMS93}
J.~J.~L. Vel{\'a}zquez.
\newblock Classification of singularities for blowing up solutions in higher
  dimensions.
\newblock {\em Trans. Amer. Math. Soc.}, 338(1):441--464, 1993.

\bibitem[Vel93b]{VINDIANA93}
J.~J.~L. Vel{\'a}zquez.
\newblock Estimates on the $(n-1)$-dimensional {H}ausdorff measure of the
  blow-up set for a semilinear heat equation.
\newblock {\em Indiana Univ. Math. J.}, 42(2):445--476, 1993.

\bibitem[Zaa98]{ZAIHPANL98}
H.~Zaag.
\newblock Blow-up results for vector-valued nonlinear heat equations with no
  gradient structure.
\newblock {\em Ann. Inst. H. Poincar\'e Anal. Non Lin\'eaire}, 15(5):581--622,
  1998.

\bibitem[Zaa02a]{ZIHP02}
H.~Zaag.
\newblock On the regularity of the blow-up set for semilinear heat equations.
\newblock {\em Ann. Inst. H. Poincar\'e Anal. Non Lin\'eaire}, 19(5):505--542,
  2002.

\bibitem[Zaa02b]{ZCMP02}
H.~Zaag.
\newblock One-dimensional behavior of singular {$N$}-dimensional solutions of
  semilinear heat equations.
\newblock {\em Comm. Math. Phys.}, 225(3):523--549, 2002.

\bibitem[Zaa02c]{ZMME02}
H.~Zaag.
\newblock Regularity of the blow-up set and singular behavior for semilinear
  heat equations.
\newblock In {\em Mathematics \& mathematics education (Bethlehem, 2000)},
  pages 337--347. World Sci. Publishing, River Edge, NJ, 2002.

\end{thebibliography}
\bibliographystyle{alpha}

\end{document}